\documentclass{amsart}
\usepackage{amssymb}

\theoremstyle{definition}

\theoremstyle{remark}

\numberwithin{equation}{section}

\input{tcilatex}

\begin{document}
\title{ General H\"{o}rmander and Mikhlin conditions for multipliers of
Besov spaces }
\author{Rishad Shahmurov}
\address{Department of Mathematics, Yeditepe University, Kayishdagi Caddesi,
34755 Kayishdagi, Istanbul, Turkey}
\email{shahmurov@hotmail.com}
\date{}

\begin{abstract}
Here a new condition for the geometry of Banach spaces is introduced and the
operator--valued Fourier multiplier theorems in weighted Besov spaces are
obtained. Particularly, connections between the geometry of Banach spaces
and H\"{o}rmander-Mikhlin conditions are established. As an application of
main results the regularity properties of degenerate elliptic differential
operator equations are investigated.
\end{abstract}

\maketitle




\section*{1. Introduction, notations and background}

In recent years, Fourier multiplier theorems in vector--valued function
spaces have found many applications in the theory of differential operator
equations, especially in maximal regularity of parabolic and elliptic
differential--operator equations and embedding theorems of abstract function
spaces. Operator--valued multiplier theorems in Banach--valued function
spaces have been discussed extensively in [3,8,12, 15,17 and 19]. Boundary
value problems (BVPs) for differential--operator equations (DOEs) in $H$%
--valued (Hilbert valued space) function spaces have been studied in
[1,2,6,7,9,13,14], and the references therein.

$D(\Omega ;E)$ will denote the collection of infinitely differentiable $E$%
--valued functions with compact support on $\Omega .$ Moreover, we denote a
bounded and uniformly continuous function spaces with traditional notation $%
BUC^{\theta }$ where 
\begin{equation*}
\left\Vert f\right\Vert _{BUC^{\theta }(\Omega ;E)}=\sup_{s\in \Omega
}\left\Vert f(s)\right\Vert +\sup_{\substack{ t,s\in \Omega  \\ s<t}}\frac{%
\left\Vert f(t)-f(s)\right\Vert _{E}}{\left\vert t-s\right\vert ^{\theta }}%
\text{ for }0<\theta <1
\end{equation*}%
and 
\begin{equation*}
\left\Vert f\right\Vert _{BUC^{m+\theta }(\Omega ;E)}=\sup_{s\in \Omega
}\sum\limits_{k=0}^{m}\left\Vert f^{(k)}(s)\right\Vert _{E}+\sup_{\substack{ %
t,s\in \Omega  \\ s<t}}\frac{\left\Vert f^{(m)}(t)-f^{(m)}(s)\right\Vert _{E}%
}{\left\vert t-s\right\vert ^{\theta }}.
\end{equation*}

Let $S(R^{n};E)$ denote the Schwartz class, i.e., a space of $E$--valued
rapidly decreasing smooth functions on $R^{n}.$ $S^{\dagger }(R^{n};E)$
denotes the space of continuous linear operators $L:S\rightarrow E$ equipped
with the bounded convergence topology sometimes called $E$--valued tempered
distributions. Let $\alpha =(\alpha _{1},\alpha _{2},\cdots ,\alpha _{n}),$
where $\alpha _{i}$ are integers. An $E$--valued generalized function $%
D^{\alpha }f$ is called a generalized derivative in the sense of Schwartz
distributions, if the equality 
\begin{equation*}
<D^{\alpha }f,\varphi >~=~(-1)^{|\alpha |}<f,D^{\alpha }\varphi >
\end{equation*}%
holds for all $\varphi \in S.$ It is known that 
\begin{equation*}
F(D_{x}^{\alpha }f)~=~(i\xi _{1})^{\alpha _{1}}\cdots (i\xi _{n})^{\alpha
_{n}}\hat{f},~~~D_{\xi }^{\alpha }(F(f))~=~F[(-ix_{n})^{\alpha _{1}}\cdots
(-ix_{n})^{\alpha _{n}}f]
\end{equation*}%
for all $f\in S^{\dagger }(R^{n};E).$

Let $\mathbf{C}$ be the set of complex numbers and 
\begin{equation*}
S_{\varphi }~=~\{\lambda ;~\lambda \in \mathbf{C},~|\arg \lambda |\leq
\varphi \}\cup \{0\},~~~0\leq \varphi <\pi .
\end{equation*}%
A linear operator $A$ is said to be a $\varphi $--positive in a Banach space 
$E$, if $D(A)$ is dense in $E$, and 
\begin{equation*}
\left\Vert (A+\lambda I)^{-1}\right\Vert _{B(E)}~\leq ~M(1+|\lambda |)^{-1}
\end{equation*}%
with $M>0,~\lambda \in S_{\varphi },~\varphi \in \lbrack 0,\pi )$; here $I$
is the identity operator in $E,~B(E)$ is the space of all bounded linear
operators in $E.$ Sometimes instead of $A+\lambda I$, we will write $%
A+\lambda $ and denote it by $A_{\lambda }.$

Let $E$ be a Banach space and $\gamma =\gamma (x),~x=(x_{1},x_{2},\cdots
,x_{n})\in \Omega \subset R^{n}.$ $L_{p,\gamma }(\Omega ;E)$ denotes the
space of all strongly measurable $E$--valued functions that are defined on
the measurable subset $\Omega \subset R^{n}$ with the norm 
\begin{equation*}
\begin{array}{lll}
\Vert f\Vert _{L_{p,\gamma }(\Omega ;E)} & = & \displaystyle\left( \int
\Vert f(x)\Vert _{E}^{p}\gamma (x)dx\right) ^{\frac{1}{p}},~~~1\leq p<\infty
, \\ 
\vspace{-3mm} &  &  \\ 
\Vert f\Vert _{L_{\infty ,\gamma }(\Omega ;E)} & = & \displaystyle%
\mbox{ess\,sup}_{x\in \Omega }[\Vert f(x)\Vert _{E}\gamma (x)].%
\end{array}%
\end{equation*}%
For $\gamma (x)\equiv 1$, we denote $L_{p,\gamma }(\Omega ;E)$ by $%
L_{p}(\Omega ;E).$ Note that dual of the space $L_{p,\gamma }(\Omega ;E)$ is
given by $L_{p^{\prime },\gamma ^{-1}}(\Omega ;E^{\prime })$ where $\frac{1}{%
p}+\frac{1}{p^{\prime }}=1$ and $\gamma ^{-1}(x)=\frac{1}{\gamma (x)}.$

We shall use Fourier analytic definition of weighted Besov spaces in this
study. Therefore, we need to consider some subsets $\{J_{k}\}_{k=0}^{\infty
} $ and $\{I_{k}\}_{k=0}^{\infty }$ of $R^{N}$, where 
\begin{equation*}
J_{0}~=~\left\{ t\in R^{N}:~|t|\leq 1\right\} ,~~~J_{k}~=~\left\{ t\in
R^{N}:~2^{k-1}\leq |t|\leq 2^{k}\right\} ~~~\mbox{for}~~~k\in N
\end{equation*}%
and 
\begin{equation*}
I_{0}~=~\left\{ t\in R^{N}:~|t|\leq 2\right\} ,~~~I_{k}~=~\left\{ t\in
R^{N}:~2^{k-1}\leq |t|\leq 2^{k+1}\right\} ~~~\mbox{for}~~~k\in N.
\end{equation*}%
Next, we define the unity $\{\varphi _{k}\}_{k\in N_{0}}$ of functions from $%
S(R^{N},R).$ Let $\psi \in S(R,R)$ be nonnegative function with support in $%
[2^{-1},2],$ which satisfies 
\begin{equation*}
\sum\limits_{k=-\infty }^{\infty }\psi (2^{-k}s)~=~1~~~\mbox{for}~~~s\in
R\backslash \{0\}
\end{equation*}%
and 
\begin{equation*}
\varphi _{k}(t)~=~\psi (2^{-k}|t|),~~~\varphi
_{0}(t)~=~1-\sum\limits_{k=1}^{\infty }\varphi _{k}(t)~~~\mbox{for}~~~t\in
R^{N}.
\end{equation*}%
Later, we will need the following useful properties: 
\begin{equation*}
\begin{array}{l}
\mbox{supp}~\varphi _{k}~\subset ~\bar{I}_{k}~~~\mbox{for
each}~~~k\in N_{0}, \\ 
\vspace{-3mm} \\ 
\varphi _{k}~\equiv ~0~~~\mbox{for each}~~~k<0, \\ 
\vspace{-3mm} \\ 
\displaystyle\sum\limits_{k=0}^{\infty }\varphi _{k}(s)~~~\mbox{for each}%
~~~s\in R^{N}, \\ 
\vspace{-3mm} \\ 
J_{m}\cap \mbox{supp}~\varphi _{k}~=~\emptyset ~~~\mbox{if}~~|m-k|>1, \\ 
\vspace{-3mm} \\ 
\varphi _{k-1}(s)+\varphi _{k}(s)+\varphi _{k+1}(s)~=~1~~~\mbox{for each}%
~~~s\in \mbox{supp}~\varphi _{k},~~~k\in N_{0}.%
\end{array}%
\end{equation*}

\vspace{1mm}

Let $1\leq q\leq r\leq \infty $ and $s\in R.$ The weighted Besov space is
the set of all functions $f\in S^{\prime }(R^{N},X)$ for which 
\begin{equation*}
\begin{array}{lll}
\Vert f\Vert _{B_{q,r,\gamma }^{s}(R^{N},X)}: & = & \displaystyle\left\Vert
2^{ks}\left\{ (\check{\varphi}_{k}\ast f)\right\} _{k=0}^{\infty
}\right\Vert _{l_{r}(L_{q,w_{q}}(R^{N},X))} \\ 
\vspace{-3mm} &  &  \\ 
& \equiv & \displaystyle\left\{ 
\begin{array}{ll}
\displaystyle\left[ \sum\limits_{k=0}^{\infty }2^{ksr}\Vert \check{\varphi}%
_{k}\ast f\Vert _{L_{q,\gamma }(R^{N},X)}^{r}\right] ^{\frac{1}{r}} & %
\mbox{if}~~r\neq \infty \\ 
\vspace{-3mm} &  \\ 
\displaystyle\sup\limits_{k\in N_{0}}\left[ 2^{ks}\Vert \check{\varphi}%
_{k}\ast f\Vert _{L_{q,\gamma }(R^{N},X)}\right] & \mbox{if}~~r=\infty%
\end{array}%
\right.%
\end{array}%
\end{equation*}%
is finite; here $q$ and $s$ are main and smoothness indexes respectively. It
is well known that Besov spaces has significant embedding properties. Thus
we close section with stating some of them: 
\begin{equation*}
W_{q}^{l+1}(X)\hookrightarrow B_{q,r}^{s}(X)\hookrightarrow
W_{q}^{l}(X)\hookrightarrow L_{q}(X)\text{ where }l<s<l+1,
\end{equation*}%
\begin{equation*}
B_{\infty ,1}^{s}(X)\hookrightarrow BUC^{s}(X)\hookrightarrow B_{\infty
,\infty }^{s}(X)\text{ for }s\in \mathbf{Z},
\end{equation*}%
and%
\begin{equation*}
B_{p,1}^{\frac{N}{p}}(R^{N},X)\hookrightarrow L_{\infty }(R^{N},X)\text{ for 
}s\in \mathbf{Z}.
\end{equation*}

For more detailed information see [2] and [3]. Let $E$ and $E_{0}$ be Banach
spaces so that $E_{0}$ is continuously and densely embedded in $E.$ We
define Besov-Lions spaces as follows: 
\begin{equation*}
\begin{array}{l}
\displaystyle B_{p,q}^{\left[ l\right] ,s}(R;E_{0},E)~=~\left\{ u:u\in
B_{p,q}^{s}(R;E_{0}),~~D^{\left[ l\right] }u\in B_{p,q}^{s}(R;E)\right\} ,
\\ 
\vspace{-3mm} \\ 
\displaystyle\Vert u\Vert _{B_{p,q}^{\left[ l\right] ,s}(R;E_{0,}E)}~=~\Vert
u\Vert _{B_{p,q}^{s}(R;E_{0})}+\left\Vert D^{\left[ l\right] }u\right\Vert
_{B_{p,q}^{s}(R,E)}~<~\infty .%
\end{array}%
\end{equation*}

\vspace{1mm}We will use this function spaces in embedding theorems and in
the study of degenerate elliptic equations.

\section*{2. Fourier multipliers}

In this section, we shall extend the work of Girardi and Weis [10] which
includes many classical multiplier conditions such as Mikhlin and H\"{o}%
rmander. This section is organized in a similar format as [10]. Some new
definitions and lemmas will be introduced. In this section $X$ and $Y$ are
Banach spaces over the field $C$ and $X^{\ast }$ is the dual space of $X.$
The space $B(X,Y)$ of bounded linear operators from $X$ to $Y$ is endowed
with the usual uniform operator topology. $N_{0}$ is the set of natural
numbers containing zero.

It is well known that Fourier transform $F:S(X)\rightarrow S(X)$ is defined
by 
\begin{equation*}
(Ff)(t)~\equiv~\hat{f}(t)~=~\int\limits_{R^N}\exp(-its)f(s)ds
\end{equation*}
is an isomorphism whose inverse is given by 
\begin{equation*}
(F^{-1}f)(t)~\equiv~\check{f}(t)~=~(2\pi)^{-N}\int\limits_{R^N}%
\exp(its)f(s)ds,
\end{equation*}
where $f\in S(X)$ and $t\in R^{N}.$

All the basic properties of $F$ and $F^{-1}$ that hold in the scalar--valued
case also hold in vector--valued case; however, the Housdorff--Young
inequality need not hold. Therefore, we need to define similar class of
Banach spaces that was introduced by Peetre [11].

\vspace{3mm}

\textbf{Definition 2.1.} Let $X$ be a Banach space and $1\leq p\leq 2.$ We
say $X$ has Fourier $\gamma $-type $p$ if 
\begin{equation*}
\Vert Ff\Vert _{L_{p^{\prime },\gamma ^{-1}}(R^{N},X)}~\leq ~C\Vert f\Vert
_{L_{p,\gamma }(R^{N},X)}~~~\mbox{for each}~~~f\in S(R^{N},X),
\end{equation*}%
where $\frac{1}{p}+\frac{1}{p^{\prime }}=1,~F_{p,N}(X)$ is the smallest $%
C\in \lbrack 0,\infty ]$ and $X$ has Fourier type $p$ if $\gamma =1.$

\vspace{3mm}

\textbf{Proposition 2.2.} Let $X$ be a Banach space with Fourier $\gamma $%
-type $p\in \lbrack 1,2]$ and $p\leq q\leq p^{\prime }.$ Then $X^{\ast }$
and $L_{q,\gamma }(R^{N},X)$ also have Fourier $\gamma $-type $p$ provided
both are with the same constant $F_{p,N}(X).$

\textbf{Proof.} It follows directly from the proof of [10, Proposition 2.3],
and Fourier $\gamma $-type property of $X.$~~~~\hspace{3mm}%
\hbox{\vrule
height7pt width5pt}

\vspace{3mm}

\textbf{Proposition 2.3.} Let $X$ be a Banach space with Fourier type $p\in
\lbrack 1,2]$. If $\gamma ^{-1}\in L_{\infty }(R^{N})$ then $X$ has a
Fourier $\gamma $-type $p.$

Since $\gamma ^{-1}\in L_{\infty }(R^{N})$ then there exist $C>0$ such that $%
1\leq C\gamma (t)$ a.e$.$ Hence we get 
\begin{equation*}
\begin{array}{lll}
\Vert \hat{f}\Vert _{L_{p^{\prime },\gamma ^{-1}}(R^{N},X)} & \leq  & %
\displaystyle C\left[ \int\limits_{R^{N}}\Vert \hat{f}(t)\Vert
_{X}^{p^{\prime }}dt\right] ^{\frac{1}{p^{\prime }}}\leq C^{2}\Vert f\Vert
_{L_{p,\gamma }(R^{N},X)}.%
\end{array}%
\hbox{\vrule height7pt width5pt}
\end{equation*}%
It is also possible to show that if $\gamma ^{1-\frac{1}{p}}\in L_{2}$ or $%
\gamma ^{1-\frac{1}{p}}\in L_{1}$ and $F^{-1}(\gamma ^{1-\frac{1}{p}})\in
L_{1}$ then Proposition 2.3 is still valid.  

\bigskip 

\textbf{Proposition 2.4.} Let $X$ be a Banach space, and $1\leq p<q$. If 
\begin{equation*}
\int\limits_{\Omega }\left[ \frac{\left( \tilde{\gamma}(t)\right) ^{q}}{%
\left( \gamma (t)\right) ^{p}}\right] ^{\frac{1}{q-p}}dt<\infty \eqno(1)
\end{equation*}%
then $L_{q,\gamma }(\Omega ,X)\hookrightarrow L_{p,\tilde{\gamma}}(\Omega
,X).$

\textbf{Proof.} We recall that embedding theorem for $L_{p}$ spaces is
applicable only in bounded domains; however, in weighted case embedding
works even in $R^{N}$ by imposing some conditions on weight. Suppose $f\in
L_{q,\gamma }(\Omega ,X).$ Then applying generalized H\"{o}lder inequality
and (1), we complete the proof: 
\begin{equation*}
\begin{array}{lll}
\Vert f\Vert _{L_{p,\tilde{\gamma}}(\Omega ,X)} & = & \displaystyle\left[
\int\limits_{\Omega }\Vert f(t)\left( \tilde{\gamma}(t)\right) ^{\frac{1}{p}%
}\Vert _{X}^{p}dt\right] ^{\frac{1}{p}} \\ 
\vspace{-3mm} &  &  \\ 
& \leq & \displaystyle\Vert f\Vert _{L_{q,\gamma }(\Omega ,X)}\left(
\int\limits_{\Omega }\left[ \frac{\left( \tilde{\gamma}(t)\right) ^{q}}{%
\left( \gamma (t)\right) ^{p}}\right] ^{\frac{1}{q-p}}dt\right) ^{\frac{q-p}{%
pq}}\leq C\Vert f\Vert _{L_{q,\gamma }(\Omega ,X)}.\text{ \ \ }%
\hbox{\vrule
height7pt width5pt}%
\end{array}%
\end{equation*}

The following Fourier embedding theorem plays a key role in the proof of
multiplier theorem.

\vspace{3mm}

\textbf{Theorem 2.5.} Let $X$ be a Banach space with the Fourier $\gamma $%
-type $p\in \lbrack 1,2].$ Let $1\leq q<p^{\prime },~1\leq r\leq \infty $
and $s\geq \frac{N}{u}$ where $\frac{1}{u}=\left( \frac{1}{q}-\frac{1}{%
p^{\prime }}\right) .$ Assume for each bounded domain $\Omega \subset R^{N}$ 
\begin{equation*}
\int\limits_{\Omega }\left[ \left( \gamma (t)\right) ^{q}\left( \tilde{\gamma%
}(t)\right) ^{p^{\prime }}\right] ^{\frac{1}{p^{\prime }-q}}dt<\infty .\eqno%
(2)
\end{equation*}%
Then there exists a constant $C$ depending only on $F_{p,N}(X),$ so that if $%
f\in B_{p,r,\gamma }^{s}(R^{N},X),$ 
\begin{equation*}
\left\Vert \left\{ \hat{f}.\chi _{J_{m}}\right\} _{k=0}^{\infty }\right\Vert
_{l_{r}(L_{q,\tilde{\gamma}}(R^{N},X))}~\leq ~C\Vert f\parallel
_{B_{p,r,\gamma }^{s}(R^{N},X)}.
\end{equation*}

\vspace{1mm}

Note that Theorem 2.5 remains valid if Fourier transform is replaced by the
inverse Fourier transform.

\textbf{Proof.} Let $f$ be in $B_{p,r,\gamma }^{s}(R^{N},X).$ Then for all $%
k\in N_{0},$ since $\check{\varphi}_{k}\ast f\in L_{p,\gamma }(R^{N},X)$ and 
$X$ has Fourier $\gamma $-type $p,~\varphi _{k}\cdot \hat{f}\in L_{p^{\prime
},\gamma ^{-1}}(R^{N},X).$ By using Proposition 2.4 and (2), we have 
\begin{equation*}
\hat{f}\cdot \chi _{J_{m}}~=~\left( \sum\limits_{k=m-1}^{m+1}\varphi _{k}%
\hat{f}\right) \chi _{J_{m}}\in L_{q,\tilde{\gamma}}(R^{N},X).
\end{equation*}%
If there exists a constant $C_{1}$ so that 
\begin{equation*}
\Vert \hat{f}.\chi _{J_{m}}\Vert _{L_{q,\tilde{\gamma}}(R^{N},X)}~\leq
C_{1}\sum\limits_{k=m-1}^{m+1}2^{ks}\left\Vert \hat{f}\cdot \varphi
_{k}\right\Vert _{L_{p^{\prime },\gamma ^{-1}}(R^{N},X)}\eqno(3)
\end{equation*}%
for each $m\in N_{0}$ then 
\begin{equation*}
\begin{array}{lll}
\Vert \hat{f}.\chi _{J_{m}}\Vert _{L_{q,\tilde{\gamma}}(R^{N},X)} & \leq & %
\displaystyle C_{1}\sum\limits_{k=m-1}^{m+1}2^{ks}\Vert F(\check{\varphi}%
_{k}\ast f\Vert _{L_{p^{\prime },\gamma ^{-1}}(R^{N},X)} \\ 
\vspace{-3mm} &  &  \\ 
& \leq & \displaystyle C_{1}F_{p,N}(X)\sum\limits_{k=m-1}^{m+1}2^{ks}\Vert 
\check{\varphi}_{k}\ast f\Vert _{L_{p,\gamma }(R^{N},X)}%
\end{array}%
\end{equation*}%
and so 
\begin{equation*}
\left\Vert \left\{ \hat{f}.\chi _{J_{m}}\right\} _{m=0}^{\infty }\right\Vert
_{l_{r}(L_{q,\tilde{\gamma}}(R^{N},X))}~\leq ~CF_{p,N}(X)\Vert f\Vert
_{B_{p,r,\gamma }^{s}}.
\end{equation*}%
It remains to show that (3) holds for some constant $C_{1}.$ Taking into
consideration (2) and applying generalized H\"{o}lder's inequality for each $%
m\in N_{0},$ we complete the proof: 
\begin{equation*}
\begin{array}{lll}
\Vert \hat{f}.\chi _{J_{m}}\Vert _{L_{q,\tilde{\gamma}}(X)} & \leq & %
\displaystyle\sum\limits_{k=m-1}^{m+1}\left\Vert \hat{f}\cdot \varphi
_{k}\cdot \chi _{J_{m}}\right\Vert _{L_{q,\tilde{\gamma}}(X)} \\ 
\vspace{-3mm} &  &  \\ 
& \leq & \displaystyle\sum\limits_{k=m-1}^{m+1}\left\Vert \hat{f}\varphi _{k}%
\left[ \frac{1+|\cdot |}{4}\right] ^{\frac{N}{u}}\chi _{J_{m}}\right\Vert
_{L_{p^{\prime },\gamma ^{-1}}(X)} \\ 
&  &  \\ 
\vspace{-3mm} & \times & \displaystyle\left\Vert \left[ \frac{1+|\cdot |}{4}%
\right] ^{\frac{-N}{u}}\chi _{J_{m}}\left( \gamma (\cdot )\right) ^{\frac{1}{%
p^{\prime }}}\left( \tilde{\gamma}(\cdot )\right) ^{\frac{1}{q}}\right\Vert
_{L_{u}(R)} \\ 
&  &  \\ 
& \leq & \displaystyle\sum\limits_{k=m-1}^{m+1}\left\Vert \left[ \frac{%
1+|\cdot |}{4}\right] ^{\frac{N}{u}}\chi _{J_{m}}\right\Vert _{L_{\infty
}(R)}\left\Vert \hat{f}\varphi _{k}\right\Vert _{L_{p^{\prime },\gamma
^{-1}}(X)} \\ 
& \times & \left[ \int\limits_{_{J_{m}}}\left[ \frac{1+|t|}{4}\right] ^{-N}%
\left[ \left( \gamma (t)\right) ^{q}\left( \tilde{\gamma}(t)\right)
^{p^{\prime }}\right] ^{\frac{1}{p^{\prime }-q}}dt\right] ^{\frac{1}{u}} \\ 
& \leq & \displaystyle C\sum\limits_{k=m-1}^{m+1}\left( 2^{m-1}\right) ^{%
\frac{N}{u}}\left\Vert \hat{f}\varphi _{k}\right\Vert _{L_{p^{\prime
},\gamma ^{-1}}(X)}~~~~~~\hspace{3mm} \\ 
& \leq & \displaystyle C\sum\limits_{k=m-1}^{m+1}2^{ks}\left\Vert \hat{f}%
\varphi _{k}\right\Vert _{L_{p^{\prime },\gamma ^{-1}}(X)}.\text{ \ \ \ }%
\hbox{\vrule height7pt
width5pt}%
\end{array}%
\end{equation*}

\vspace{1mm}

\textbf{Corollary 2.6.} Let $X$ be a Banach space with Fourier $\gamma $%
-type $p\in \lbrack 1,2].$ If (2) holds for $q=r=1$ and $r=q=p^{\prime }$
then the Fourier transform defines bounded operators 
\begin{equation*}
F:~B_{p,1,\gamma }^{N/p}(R^{N},X)~\rightarrow ~L_{1,\tilde{\gamma}}(R^{N},X)%
\eqno(4)
\end{equation*}%
\begin{equation*}
F:~B_{p,p^{\prime },\gamma ^{-1}}^{0}(R^{N},X)~\rightarrow ~L_{p^{\prime
},\gamma ^{-1}}(R^{N},X).\eqno(5)
\end{equation*}

For a bounded measurable function $m:R^{N}\rightarrow B(X,Y),$ its
corresponding Fourier multiplier operator $T_{m}$ is defined as follows 
\begin{equation*}
T_{m}(f)~=~F^{-1}[m(\cdot )(Ff)(\cdot )].
\end{equation*}%
In this section, we identify conditions on $m$, extending those of [10],
that 
\begin{equation*}
\Vert T_{0}f\Vert _{B_{q,r,\tilde{\gamma}}^{s}}~\leq ~C\Vert f\Vert _{B_{q,r,%
\tilde{\gamma}}^{s}}~~~\mbox{for each}~~~f\in S(X).
\end{equation*}

\vspace{1mm}

\textbf{Definition 2.7.} Let $(E(R^{N},Z),E^{\ast }(R^{N},Z^{\ast }))$ be
one of the following dual systems, where $1\leq q,~r\leq \infty $ and $s\in
R $ 
\begin{equation*}
(L_{q,\tilde{\gamma}}(Z),L_{q^{\prime },\tilde{\gamma}^{-1}}(Z^{\ast }))~~~%
\mbox{or}~~~(B_{q,r,\tilde{\gamma}}^{s}(Z),B_{q^{\prime },r^{\prime },\tilde{%
\gamma}^{-1}}^{-s}(R^{N},Z^{\ast })).
\end{equation*}%
A bounded measurable function $m:R^{N}\rightarrow B(X,Y)$ is called a
Fourier multiplier from $E(X)$ to $E(Y)$ if there is a bounded linear
operator 
\begin{equation*}
T_{m}:~E(X)\rightarrow E(Y)
\end{equation*}%
such that 
\begin{equation*}
T_{m}(f)~=~F^{-1}[m(\cdot )(Ff)(\cdot )]~~~\mbox{for each}~~~f\in S(X),\eqno%
(6)
\end{equation*}%
\begin{equation*}
T_{m}~~~\mbox{is}~~~\sigma (E(X),E^{\ast }(X^{\ast }))~~~\mbox{to}~~~\sigma
(E(Y),E^{\ast }(Y^{\ast }))~~~\mbox{continuous.}\eqno(7)
\end{equation*}%
The uniquely determined operator $T_{m}$ is the Fourier multiplier operator
induced by $m.$

\vspace{3mm}

\textbf{Remark 2.8.} If $T_{m}\in B(E(X),E(Y))$ and $T_{m}^{\ast }$ maps $%
E^{\ast }(Y^{\ast })$ into $E^{\ast }(X^{\ast })$ then $T_{m}$ satisfies the
continuity condition (7).

\vspace{3mm}

\textbf{Lemma 2.9.} Suppose $k\in L_{1,\tilde{\gamma}}(R^{N},B(X,Y))$ and
the weight function $\tilde{\gamma}$ satisfies the following condition:%
\begin{equation*}
\sup_{t\in R^{N}}\frac{\tilde{\gamma}(t)}{\tilde{\gamma}(t-s)}\leq C_{1}%
\tilde{\gamma}(s)\text{ for all }s\in R^{N}.
\end{equation*}%
Assume there exist $C_{2}$ so that 
\begin{equation*}
\Vert k(\cdot )x\Vert _{L_{1,\tilde{\gamma}}(Y)}\leq ~C_{2}\Vert x\Vert _{X}%
\text{ for all }x\in X\text{ }\eqno(8)
\end{equation*}%
and $C_{3}$ so that 
\begin{equation*}
\Vert k^{\ast }y^{\ast }\Vert _{L_{1,\tilde{\gamma}}(X^{\ast })}~\leq
~C_{3}\Vert y^{\ast }\Vert _{Y^{\ast }}\text{ for all }y^{\ast }\in Y^{\ast
}.\eqno(9)
\end{equation*}%
Then for $1\leq q\leq \infty $ the convolution operator 
\begin{equation*}
K:~L_{q,\tilde{\gamma}}(R^{N},X)~\rightarrow ~L_{q,\tilde{\gamma}}(R^{N},Y)
\end{equation*}%
defined by 
\begin{equation*}
(Kf)(t)~=~\int\limits_{R^{N}}k(t-s)f(s)ds~~~\mbox{for}~~~t\in R^{N}
\end{equation*}%
satisfies $\Vert K\Vert _{L_{q,\tilde{\gamma}}\rightarrow L_{q,\tilde{\gamma}%
}}\leq C_{1}C_{2}^{1-\frac{1}{q}}C_{3}^{\frac{1}{q}}.$

\textbf{Proof}. We shall prove this lemma in a similar manner as [10, Lemma
4.5], applying vector-valued Stein-Weiss interpolation theorem $[16,\S %
1.18.5]$ instead of [12, Theorem 5.1.2]. By using (8), we have the assertion
for $q=1.$ Really 
\begin{equation*}
\begin{array}{lll}
\Vert (Kf)(t)\Vert _{L_{1,\tilde{\gamma}}(Y)} & \leq & \displaystyle%
\int\limits_{R^{N}}\int\limits_{R^{N}}\Vert k(t-s)f(s)\Vert _{Y}\tilde{\gamma%
}(t)dt\,ds \\ 
\vspace{-3mm} &  &  \\ 
& \leq & \displaystyle\int\limits_{R^{N}}\left[ \int\limits_{R^{N}}\Vert
k(t-s)f(s)\Vert _{Y}\tilde{\gamma}(t-s)dt\right] \sup_{t\in R^{N}}\frac{%
\tilde{\gamma}(t)}{\tilde{\gamma}(t-s)}\,ds \\ 
&  & ~ \\ 
& \leq & \displaystyle C_{1}C_{2}\int\limits_{R^{N}}\Vert f(s)\Vert _{X}%
\tilde{\gamma}(s)\,ds\leq C_{1}C_{2}\Vert f(s)\Vert _{L_{1,\tilde{\gamma}%
}(X)}.%
\end{array}%
\end{equation*}%
If $f\in L_{\infty ,\tilde{\gamma}}(Y),~y^{\ast }\in Y^{\ast }$ and $t\in
R^{N}$ then $\left\Vert K\right\Vert _{L_{\infty ,\tilde{\gamma}}\rightarrow
L_{\infty ,\tilde{\gamma}}}\leq C_{1}C_{3}:$%
\begin{equation*}
\begin{array}{lll}
|<y^{\ast },(Kf)(t)\tilde{\gamma}(t)>_{Y}| & \leq & \displaystyle%
\int\limits_{R^{N}}|<k(t-s)^{\ast }y^{\ast }\tilde{\gamma}(t),f(s)>_{X}|ds
\\ 
\vspace{-3mm} &  &  \\ 
& \leq & \displaystyle\int\limits_{R^{N}}\Vert k(t-s)^{\ast }y^{\ast }\Vert
_{X^{\ast }}\tilde{\gamma}(t-s)\frac{\tilde{\gamma}(t)}{\tilde{\gamma}(t-s)}%
\Vert f(s)\Vert _{X}ds \\ 
\vspace{-3mm} &  &  \\ 
& \leq & \displaystyle~C_{1}C_{3}\Vert y^{\ast }\Vert _{Y^{\ast }}\Vert
f\Vert _{L_{\infty ,\tilde{\gamma}}(X)}.%
\end{array}%
\end{equation*}%
In view of $[16,\S 1.18.5]$, we conclude that $\Vert K\Vert _{L_{_{q,\tilde{%
\gamma}}}\rightarrow L_{q,\tilde{\gamma}}}\leq C_{1}C_{2}^{1-\frac{1}{q}%
}C_{3}^{\frac{1}{q}}$ for $1\leq q\leq \infty .$~~~~\hspace{3mm}%
\hbox{\vrule
height7pt width5pt}

\vspace{3mm}

\textbf{Proposition 2.10.} Let $E$ be a Banach space, $1\leq p<\infty $ and $%
\gamma $ be a positive measurable function on an open subset $\Omega $ of $%
R^{n},$ and essentially bounded on a compact subsets of $\Omega .$ Then $%
D(\Omega ;E)\hookrightarrow L_{p,\gamma }(\Omega ;E).$

\textbf{Proof.} For $u\in L_{p,\gamma }(\Omega ;E)$ and $n\in \mathbf{N}$
let $u_{n}:\Omega \rightarrow E$ be such that 
\begin{equation*}
u_{n}=\left\{ 
\begin{array}{ll}
u(x) & \mbox{if}~~~\Vert u(x)\Vert \leq n \\ 
0 & \mbox{if}~~~\Vert u(x)\Vert >n.%
\end{array}%
\right.
\end{equation*}

By the dominated convergence theorem $\lim\limits_{n\rightarrow \infty
}\Vert u-u_{n}\Vert _{L_{p,\gamma }\left( \Omega ;E\right) }=0,$ hence a
compactly supported function can be approximated by bounded compactly
supported functions belonging to $L_{p}(\Omega ;E).$ From the proof of the
denseness theorem (classical case), it follows that if $u$ is a compactly
supported function belonging to $L_{p}(\Omega ;E)$ then there exists a
compact subset $K\subset \Omega ,$ with $\mbox{supp}\,u\subseteq K,$ and a
sequence of functions $u_{n}\in D(\Omega ;E),$ with $\mbox{supp}%
~u_{n}\subseteq K,$ such that $\lim_{n\rightarrow \infty }\Vert u-u_{n}\Vert
_{L_{p}(\Omega ;E)}=0;$ since 
\begin{equation*}
\Vert u-u_{n}\Vert _{L_{p,\gamma }(\Omega ;E)}~=~\left( \int\limits_{K}\Vert
u(x)-u_{n}(x)\Vert ^{p}\gamma (x)dx\right) ^{\frac{1}{p}}~\leq ~\left(
\sup\limits_{x\in K}\gamma (x)\right) ^{\frac{1}{p}}\Vert u-u_{n}\Vert
_{L_{p}(\Omega ;E)}
\end{equation*}%
we have $\lim\limits_{n\rightarrow \infty }\Vert u-u_{n}\Vert _{L_{p,\gamma
}(\Omega ;E)}=0.$~~~~\hspace{3mm}\hbox{\vrule
height7pt width5pt}

\vspace{3mm}

\textbf{Condition 1.} Let $p^{\prime }$ be a dual pair of $p$ (Fourier $%
\gamma $-type of a Banach spaces $X$ and $Y$). Suppose $\gamma $ is
measurable on each open subset $\Omega \subset R^{N},$ essentially bounded
on each compact subset $\Omega \subset R^{N}~$and for each $t\in R^{N}:$%
\begin{equation*}
(i)\sup_{t\in R^{N}}\frac{\tilde{\gamma}(t)}{\tilde{\gamma}(t-s)}\leq C%
\tilde{\gamma}(s)\text{ for all }s\in R^{N};\eqno(10)
\end{equation*}%
\begin{equation*}
(ii)\int\limits_{\Omega }\left[ \left( \gamma (t)\right) ^{1-\frac{1}{p}}%
\tilde{\gamma}(t)\right] ^{p}dt<\infty \text{ for each }\Omega \subset R^{N},%
\text{ vol}\left( \Omega \right) <\infty .
\end{equation*}

\vspace{3mm}

\textbf{Example 1. (a.) }It is easy to see that exponential functions
satisfy Condition 1.

\textbf{(b.) }As a second example we can give e.g. polynomials functions of
the form $(1+\left\vert x\right\vert )^{k}.$

In $\left[ 18\right] $ author studied FMT in $L_{p,\gamma }(R^{N},l_{p})$
for $p\in \left( 1,\infty \right) .$ Particularly, it was shown that
choosing weight functions in the following form%
\begin{eqnarray*}
(i)\gamma &=&\left\vert x\right\vert ^{\alpha }\text{, }-1<\alpha <p-1\text{%
, } \\
(ii)\gamma &=&\prod\limits_{k=1}^{N}\left( 1+\sum\limits_{j=1}^{n}\left\vert
x_{j}\right\vert ^{\alpha _{jk}}\right) ^{\beta _{k}}\text{, }\alpha
_{jk}\geq 0\text{, }N\in \mathbf{N}\text{, }\beta _{k}\in R
\end{eqnarray*}%
it is possible to establish boundedness of Fourier multiplier operator.

\vspace{3mm}

\textbf{Theorem 2.13.} Let $X$ and $Y$ be Banach spaces with Fourier $\gamma 
$-type $p\in \lbrack 1,2]$. Assume Condition 1 holds. Then there is a
constant $C$ depending only on $F_{p,N}(X)$ and $F_{p,N}(Y)$ so that if 
\begin{equation*}
m~\in ~B_{p,1,\gamma }^{\frac{N}{p}}(R^{N},B(X,Y))
\end{equation*}%
then $m$ is a Fourier multiplier from $L_{q,\tilde{\gamma}}(R^{N},X)$ to $%
L_{q,\tilde{\gamma}}(R^{N},Y)$ with 
\begin{equation*}
\Vert T_{m}\Vert _{L_{q,\tilde{\gamma}}(R^{N},X)\rightarrow L_{q,\tilde{%
\gamma}}(R^{N},Y)}~\leq ~CM_{p}(m)~~~\mbox{for each}~~~q\in \lbrack 1,\infty
]\eqno(11)
\end{equation*}%
where 
\begin{equation*}
M_{p,\gamma }(m)~=~\inf \left\{ \Vert m(a\cdot )\Vert _{B_{p,1,\gamma }^{%
\frac{N}{p}}(R^{N},B(X,Y))}:~a>0\right\} .
\end{equation*}

\textbf{Proof. }The key points in this proof are the fact $\left( 4\right) $
and Lemma 2.9. As in the proof of $[10,$Theorem $4.3]$ we assume in addition
that $m\in S\left( B\left( X,Y\right) \right) .$ Hence, $\check{m}\in
S\left( B\left( X,Y\right) \right) .$ Since, $F^{-1}\left[ m\left( a\cdot
\right) x\right] \left( s\right) =a^{-N}\check{m}\left( \frac{s}{a}\right) x$%
~choosing an appropriate $a$ and using $\left( 1\right) $ we obtain 
\begin{equation*}
\begin{array}{lll}
\vspace{-3mm} &  & \displaystyle\left\Vert \check{m}\left( \cdot \right)
x\right\Vert _{L_{1,\tilde{\gamma}}\left( Y\right) }=\left\Vert \left[
m\left( a\cdot \right) x\right] ^{\vee }\right\Vert _{L_{1,\tilde{\gamma}%
}\left( Y\right) } \\ 
&  &  \\ 
& \leq & \displaystyle C_{1}\left\Vert m\left( a\cdot \right) \right\Vert
_{B_{p,1,\gamma }^{\frac{N}{p}}}\left\Vert x\right\Vert _{X}\leq
2C_{1}M_{p,\gamma }(m)\left\Vert x\right\Vert _{X}%
\end{array}%
\end{equation*}%
~

where $C_{1}$ depends only on $F_{p,N}(Y).$~If $m\in S\left( B\left(
X,Y\right) \right) $ then $\left[ m(\cdot )^{\ast }\right] ^{\vee }=\left[ 
\check{m}(\cdot )\right] ^{\ast }\in S\left( B\left( Y^{\ast },X^{\ast
}\right) \right) $ and $M_{p,\gamma }(m)=M_{p,\gamma }(m^{\ast }).$ Thus, in
a similar manner as above, we have$\ \ \ \ \ \ \ \ \ \ \ \ \ \ \ \ \ \ \ \ \
\ \ \ \ \ \ \ \ \ \ \ \ \ \ \ \ \ \ \ \ \ \ \ \ \ \ \ $%
\begin{equation*}
\ 
\begin{array}{lll}
\vspace{-3mm} &  & \displaystyle\left\Vert \left[ \check{m}\left( \cdot
\right) \right] ^{\ast }y^{\ast }\right\Vert _{L_{1,\tilde{\gamma}}\left(
Y\right) }\leq 2C_{2}M_{p,\gamma }(m)\left\Vert y^{\ast }\right\Vert
_{Y^{\ast }}%
\end{array}%
\end{equation*}
for some constant $C_{2}$ depends on $F_{p,N}(X^{\ast }).$ Since, we have%
\begin{equation*}
\left\Vert \check{m}\left( \cdot \right) x\right\Vert _{L_{1,\tilde{\gamma}%
}\left( Y\right) }\leq 2C_{1}M_{p,\gamma }(m)\left\Vert x\right\Vert _{X}
\end{equation*}%
and%
\begin{equation*}
\left\Vert \left[ \check{m}\left( \cdot \right) \right] ^{\ast }y^{\ast
}\right\Vert _{L_{1,\tilde{\gamma}}\left( Y\right) }\leq 2C_{2}M_{p,\gamma
}(m)\left\Vert y^{\ast }\right\Vert _{Y^{\ast }}
\end{equation*}%
by Lemma 2.9 we can conclude 
\begin{equation*}
\left( T_{m}f\right) \left( t\right) )=\dint\limits_{R^{N}}\check{m}\left(
t-s\right) f\left( s\right) ds
\end{equation*}%
satisfies%
\begin{equation*}
\left\Vert T_{m}f\right\Vert _{L_{q,\tilde{\gamma}}\left( R^{N},Y\right)
}\leq CM_{p,\gamma }(m)\left\Vert f\right\Vert _{L_{q,\tilde{\gamma}}\left(
R^{N},X\right) .}
\end{equation*}%
Now, taking into account the fact that $S\hookrightarrow B_{p,1,\gamma }^{%
\frac{N}{p}}$ and using the same reasoning as in the proof of $[10,$Theorem $%
4.3]$ one can easily prove for the general case $m\in B_{p,1,\gamma }^{\frac{%
N}{p}}$ and that $T_{m}$ satisfies $\left( 7\right) .$ \ 
\hbox{\vrule
height7pt width5pt}

\vspace{3mm}

\textbf{Theorem 2.14.} Let $X$ and $Y$ be Banach spaces with Fourier $\gamma 
$-type $p\in \lbrack 1,2].$ Assume Condition 1 holds. Then there exist a
constant $C$ depending only on $F_{p,N}(X)$ and $F_{p,N}(Y)$ so that if $%
m:R^{N}\rightarrow B(X,Y)$ satisfy 
\begin{equation*}
\varphi _{k}\cdot m\in B_{p,1,\gamma }^{\frac{N}{p}}(R^{N},B(X,Y))~~~%
\mbox{and}~~~M_{p,\gamma }(\varphi _{k}\cdot m)~\leq ~A\eqno(12)
\end{equation*}%
then $m$ is Fourier multiplier from $B_{q,r,\tilde{\gamma}}^{s}(R^{N},X)$ to 
$B_{q,r,\tilde{\gamma}}^{s}(R^{N},Y)$ and $\Vert T_{m}\Vert \leq CA$ for
each $s\in R$ and $r\in \lbrack 1,\infty ].$

\textbf{Proof.} By using Theorem 2.13 we shall prove this theorem in a
similar manner as $\left[ 10,\text{ Theorem }4.8\right] $~.~Really, since $%
\varphi _{k}\cdot m\in B_{p,1,\gamma }^{\frac{N}{p}}(R^{N},B(X,Y)),$ Theorem
2.13 ensures that%
\begin{equation*}
\left\Vert T_{m\varphi _{k}}f\right\Vert _{L_{q,\tilde{\gamma}%
}(R^{N},Y)}\leq C~M_{p}(\varphi _{k}\cdot m)~\left\Vert f\right\Vert _{L_{q,%
\tilde{\gamma}}(R^{N},X)}\leq ~CA\left\Vert f\right\Vert _{L_{q,\tilde{\gamma%
}}(R^{N},X)}.
\end{equation*}%
In the introduction we defined function $\psi _{k}=\varphi _{k-1}+\varphi
_{k}+\varphi _{k+1}$ that is equal to $1$ on supp$\varphi _{k}.$ Thus, 
\begin{equation*}
\begin{array}{lll}
\vspace{-3mm} &  & \displaystyle\left\Vert T_{m\psi _{k}}f\right\Vert _{L_{q,%
\tilde{\gamma}}\left( Y\right) }\leq \left\Vert T_{m\varphi
_{k-1}}f\right\Vert _{L_{q,\tilde{\gamma}}\left( Y\right) } \\ 
&  &  \\ 
& + & \displaystyle\left\Vert T_{m\varphi _{k}}f\right\Vert _{L_{q,\tilde{%
\gamma}}\left( Y\right) }+\left\Vert T_{m\varphi _{k+1}}f\right\Vert _{L_{q,%
\tilde{\gamma}}\left( Y\right) } \\ 
&  &  \\ 
\vspace{-3mm} & \leq & \displaystyle3CA\left\Vert f\right\Vert _{L_{q,\tilde{%
\gamma}}(R^{N},X)}.%
\end{array}%
\end{equation*}%
Let $T_{0}:S(X)\rightarrow S^{\prime }(Y)$ be defined as%
\begin{equation*}
T_{0}f=F^{-1}\left[ m(\cdot )\left( Ff\right) (\cdot )\right] .
\end{equation*}%
From the proof of $[10,$ Theorem $4.3]$ we know that%
\begin{equation*}
\check{\varphi}_{k}\ast T_{0}f=T_{m\psi _{k}}\left( \check{\varphi}_{k}\ast
f\right) .
\end{equation*}%
Hence,%
\begin{equation*}
\begin{array}{lll}
\vspace{-3mm} &  & \displaystyle\left\Vert \check{\varphi}_{k}\ast
T_{0}f\right\Vert _{L_{q,\tilde{\gamma}}\left( Y\right) }=\left\Vert
T_{m\psi _{k}}\left( \check{\varphi}_{k}\ast f\right) \right\Vert _{L_{q,%
\tilde{\gamma}}\left( Y\right) } \\ 
&  &  \\ 
& \leq & \displaystyle3CA\left\Vert f\right\Vert _{L_{q,\tilde{\gamma}%
}(R^{N},X)}%
\end{array}%
\end{equation*}%
and%
\begin{equation*}
\begin{array}{lll}
\vspace{-3mm} &  & \displaystyle\dsum\limits_{k=0}^{\infty
}2^{ksr}\left\Vert \check{\varphi}_{k}\ast T_{0}f\right\Vert _{L_{q,\tilde{%
\gamma}}\left( Y\right) }^{r}\leq 3CA\dsum\limits_{k=0}^{\infty
}2^{ksr}\left\Vert f\right\Vert _{L_{q,\tilde{\gamma}}(R^{N},X)}^{r}.%
\end{array}%
\end{equation*}%
Hence, we obtain%
\begin{equation*}
\left\Vert T_{0}f\right\Vert _{B_{q,r,\tilde{\gamma}}^{s}(R^{N},Y)}\leq
3CA\left\Vert f\right\Vert _{B_{q,r,\tilde{\gamma}}^{s}(R^{N},X)}
\end{equation*}%
for $1\leq q<\infty .$ If $q,r<\infty $ then $\mathring{B}_{q,r,\tilde{\gamma%
}}^{s}=B_{q,r,\tilde{\gamma}}^{s}.$ Therefore, it remains to show the cases $%
q=\infty $ and $r=\infty $ and the weak continuity condition $\left(
7\right) .$ The case $r=\infty $ and the weak continuity condition $\left(
7\right) $ can be proved in a similar manner as $[10,$ Theorem $4.3].$%
\vspace{3mm}~~~\hbox{\vrule height7pt width5pt}

In the next section, we apply Theorem 2.14 to degenerate DOEs. However,
checking assumptions of the theorem for multiplier functions is not
practical. Therefore, we prove a lemma that makes Theorem 2.14 more
applicable.

\vspace{3mm}

\textbf{Lemma 2.15.} Let $\frac{N}{p}<l\in N,$ $u\in \lbrack p,\infty ]$ and%
\begin{equation*}
\mbox{where}~~~\frac{1}{p}=\frac{1}{u}+\frac{1}{\tilde{u}}.\eqno(13)
\end{equation*}%
Moreover, suppose $X$ and $Y$ are Banach spaces having Fourier $\gamma $%
-type $p\in \lbrack 1,2]$ and Condition 1 holds. If $m\in
C^{l}(R^{N},B(X,Y)) $ satisfies the following 
\begin{equation*}
\Vert \gamma ^{\frac{1}{p}}(\cdot )D^{\alpha }m(\cdot )|_{I_{0}}\Vert
_{L_{u}(B(X,Y))}~\leq ~A,~~~\Vert \gamma ^{\frac{1}{p}}(\cdot )D^{\alpha
}m(2^{k-1}\cdot )|_{I_{1}}\Vert _{L_{u}(B(X,Y))}~\leq ~A
\end{equation*}%
for each $\alpha \in N_{0}^{N},~|\alpha |\leq l,$ then $m$ satisfies
conditions of Theorem 2.14.

\textbf{Proof.} By using the fact that $W_{p,\gamma
}^{l}(R^{N},B(X,Y))\subset B_{p,1,\gamma }^{\frac{N}{p}}(R^{N},B(X,Y))$ for $%
\frac{N}{p}<l$ and applying Holder's inequality we get desired result~~~%
\begin{eqnarray*}
M_{p,\gamma }(\varphi _{0}\cdot m) &\leq &K\Vert \varphi _{0}m\Vert
_{W_{p,\gamma }^{l}}\leq K\sum\limits_{|\alpha |\leq l}\sum\limits_{\beta
\leq \alpha }\left\Vert \left( 
\begin{array}{c}
\alpha \\ 
\beta%
\end{array}%
\right) D^{\beta }\varphi _{0}\left( \gamma ^{\frac{1}{p}}D^{\alpha -\beta
}m\right) \right\Vert _{L_{p}} \\
&\leq &K\sum\limits_{|\alpha |\leq l}\sum\limits_{\beta \leq \alpha }\left( 
\begin{array}{c}
\alpha \\ 
\beta%
\end{array}%
\right) \left\Vert D^{\beta }\varphi _{0}\right\Vert _{L_{\tilde{u}%
}(R^{N})}\cdot \sum\limits_{|\alpha |\leq l}\sum\limits_{\beta \leq \alpha
}\left\Vert \gamma ^{\frac{1}{p}}D^{\alpha }m|_{I_{0}}\right\Vert _{L_{u}} \\
&\leq &KAC_{\varphi _{0}}
\end{eqnarray*}%
and%
\begin{eqnarray*}
M_{p,\gamma }(\varphi _{k}\cdot m) &\leq &\Vert \varphi _{k}(2^{k-1}\cdot
)m(2^{k-1}\cdot )\Vert _{B_{p,1,\gamma }^{\frac{N}{p}}}=\Vert \varphi
_{1}\left( \cdot \right) m(2^{k-1}\cdot )\Vert _{B_{p,1,\gamma }^{\frac{N}{p}%
}} \\
&\leq &K\Vert \varphi _{1}\left( \cdot \right) m(2^{k-1}\cdot )\Vert
_{W_{p,\gamma }^{l}}\leq KAC_{\varphi _{1}}.
\end{eqnarray*}

\hbox{\vrule height7pt width5pt}

We close this section with two very important corollaries that provide
different sufficient conditions for $B_{q,r,\tilde{\gamma}}^{s}$-regularity
of (6). As a matter of fact these conditions are slightly modified versions
of H\"{o}rmander and Mikhlin conditions.

\vspace{3mm}

\textbf{Corollary 2.16. (}FMT via \textit{H\"{o}rmander} condition\textbf{)}
Suppose $X$ and $Y$ have Fourier $\gamma $-type $p\in \lbrack 1,2]$ and
Condition 1 holds. If $m\in C^{l}(R^{N},B(X,Y))$ satisfies 
\begin{equation*}
\left[ \dint\limits_{\left\vert t\right\vert \leq 2}\left\Vert D^{\alpha
}m(t)\right\Vert ^{p}\gamma (t)dt\right] ^{\frac{1}{p}}\leq ~A
\end{equation*}%
and 
\begin{equation*}
\left[ R^{-N}\dint\limits_{R\leq \left\vert t\right\vert \leq 4R}\left\Vert
D^{\alpha }m(t)\right\Vert ^{p}\gamma (t)dt\right] ^{\frac{1}{p}}\leq
~AR^{-\left\vert \alpha \right\vert }
\end{equation*}%
for each multi--index $\alpha $ with $|\alpha |\leq \left\lceil \frac{N}{p}%
\right\rceil +1$ then $m$ is Fourier multiplier from $B_{q,r,\tilde{\gamma}%
}^{s}(R^{N},X)$ to $B_{q,r,\tilde{\gamma}}^{s}(R^{N},Y)$ for each $s\in R$
and $r\in \lbrack 1,\infty ].$

\textbf{Proof. }Choosing $u=p$ in the Lemma 2.15 we get assertions of
corollary. \ \hbox{\vrule height7pt width5pt}

\vspace{3mm}

\textbf{Corollary 2.17. }(FMT via Mikhlin condition)\textbf{\ }Assume $X$
and $Y$ are Banach spaces with Fourier $\gamma $-type $p\in \lbrack 1,2]$
and Condition 1 holds$.$ If $m\in C^{l}(R^{N},B(X,Y))$ satisfies 
\begin{equation*}
\left\Vert \gamma ^{\frac{1}{p}}(t)(1+|t|)^{|\alpha |}D^{\alpha
}m(t)~\right\Vert _{L_{\infty }(R^{N},B(X,Y))}\leq ~A
\end{equation*}%
for each multi--index $\alpha $ with $|\alpha |\leq l=\left\lceil \frac{N}{p}%
\right\rceil +1,$ then $m$ is Fourier multiplier from $B_{q,r,\tilde{\gamma}%
}^{s}(R^{N},X)$ to $B_{q,r,\tilde{\gamma}}^{s}(R^{N},Y)$ for each $s\in
R,~r,~q\in \lbrack 1,\infty ].$

\textbf{Proof. }Choosing $u=\infty $ in the Lemma 2.15 one can prove this
result in a similar way as [10, Corollary 4.11].~~%
\hbox{\vrule height7pt
width5pt}

The following result is a special case of Corollary 2.17. Choosing $\gamma
=1 $ we obtain a sufficient condition for the multipliers of weighted Besov
spaces.

\vspace{3mm}

\textbf{Corollary 2.18. }Assume $X$ and $Y$ are Banach spaces with Fourier
type $p$ and 
\begin{equation*}
(i)\sup_{t\in R^{N}}\frac{\tilde{\gamma}(t)}{\tilde{\gamma}(t-s)}\leq C%
\tilde{\gamma}(s)\text{ for all }s\in R^{N}
\end{equation*}%
\begin{equation*}
(ii)\int\limits_{\Omega }\left[ \tilde{\gamma}(t)\right] ^{p}dt<\infty \text{
for each }\Omega \subset R^{N},\text{ vol}\left( \Omega \right) <\infty .
\end{equation*}%
If $m\in C^{l}(R^{N},B(X,Y))$ satisfies 
\begin{equation*}
\left\Vert (1+|t|)^{|\alpha |}D^{\alpha }m(t)~\right\Vert _{L_{\infty
}(R^{N},B(X,Y))}\leq ~A
\end{equation*}%
for each multi--index $\alpha $ with $|\alpha |\leq l=\left\lceil \frac{N}{p}%
\right\rceil +1,$ then $m$ is Fourier multiplier from $B_{q,r,\tilde{\gamma}%
}^{s}(R^{N},X)$ to $B_{q,r,\tilde{\gamma}}^{s}(R^{N},Y)$ for each $s\in
R,~r,~q\in \lbrack 1,\infty ].$

\section*{3. Differential Embeddings}

In the present section, by using Corollary 2.18 we shall prove continuity of
the following embedding%
\begin{equation*}
D^{\alpha }:B_{q,r,\gamma }^{l,s}\left( R^{N};E\left( A\right) ,E\right)
\subset B_{q,r,\gamma }^{s}\left( R^{N};E\right) .
\end{equation*}%
In the next section we will apply above result to non degenerate elliptic
equations.

\vspace{3mm}

\textbf{Condition 2. }Assume a positive weight function $\gamma $ satisfies
the following: 
\begin{equation*}
(i)\sup_{t\in R^{N}}\frac{\gamma (t)}{\gamma (t-s)}\leq C\gamma (s)\text{
for all }s\in R
\end{equation*}%
\begin{equation*}
(ii)\int\limits_{\Omega }\left[ \gamma (t)\right] ^{p}dt<\infty \text{ for
each }\Omega \subset R,\text{ vol}\left( \Omega \right) <\infty .
\end{equation*}

\vspace{3mm}

\textbf{Theorem 3.1}. Suppose Condition 2 holds and $0<$\ $h\leq
h_{0}<\infty $. Let $E$ be a Banach space with Fourier type $p$ and $A$ be a 
$\varphi $-positive operator in $E,$ where $\varphi \in \left( 0\right.
,\left. \pi \right] .$ If $\alpha =\left( \alpha _{1},\alpha _{2},...,\alpha
_{N}\right) ,$ $x=\frac{\left\vert \alpha \right\vert }{l}\leq 1$ and $0<\mu
\leq 1-x$ then the following embedding 
\begin{equation*}
D^{\alpha }:B_{q,r,\gamma }^{l,s}\left( R^{N};E\left( A\right) ,E\right)
\subset B_{q,r,\gamma }^{s}\left( R^{N};E\left( A^{1-x-\mu }\right) \right)
\end{equation*}%
is continuous and there exists a positive constant\ $C$ such that%
\begin{eqnarray*}
&&\left\Vert D^{\alpha }u\right\Vert _{B_{q,r,\gamma }^{s}\left(
R^{N};E\left( A^{1-x-\mu }\right) \right) } \\
&\leq &C_{\mu }\left[ h^{\mu }\left\Vert u\right\Vert _{B_{q,r,\gamma
}^{l,s}\left( R^{N};E\left( A\right) ,E\right) }+h^{-\left( 1-\mu \right)
}\left\Vert u\right\Vert _{B_{q,r,\gamma }^{s}\left( R^{N};E\right) }\right]
\end{eqnarray*}%
for all $u\in B_{q,r,\gamma }^{l,s}\left( R^{N};E\left( A\right) ,E\right) .$

\textbf{Proof}. Since $A$ is constant and closed operator$,$ we have%
\begin{equation*}
\begin{array}{lll}
\left\Vert D^{\alpha }u\right\Vert _{B_{q,r,\gamma }^{s}\left( R^{N};E\left(
A^{1-x-\mu }\right) \right) } & = & \left\Vert A^{1-x-\mu }D^{\alpha
}u\right\Vert _{B_{q,r,\gamma }^{s}\left( R^{N};E\right) } \\ 
& \backsim & \left\Vert F^{-\shortmid }\left( i\xi \right) ^{\alpha
}A^{1-x-\mu }Fu\right\Vert _{B_{q,r,\gamma }^{s}\left( R^{N};E\right) }.%
\end{array}%
\end{equation*}%
(The symbol $\backsim $ indicates norm equivalency). In a similar manner,
from definition of $B_{q,r,\gamma }^{l,s}\left( R^{N};E_{0},E\right) $ we
have%
\begin{equation*}
\left\Vert u\right\Vert _{B_{q,r,\gamma }^{l,s}\left( R^{N};E_{0},E\right)
}\sim \left\Vert Au\right\Vert _{B_{q,r,\gamma }^{s}\left( R^{N};E\right)
}+\sum\limits_{k=1}^{N}\left\Vert F^{-1}\xi _{k}^{l}\hat{u}\right\Vert
_{B_{q,r,\gamma }^{s}}.
\end{equation*}%
By virtue of above relations, it is sufficient to prove 
\begin{eqnarray*}
&&\left\Vert F^{-\shortmid }\left[ \left( i\xi \right) ^{\alpha }A^{1-x-\mu }%
\hat{u}\right] \right\Vert _{B_{q,r,\gamma }^{s}\left( R^{N};E\right) } \\
&\leq &C\left[ \left\Vert F^{-\shortmid }A\hat{u}\right\Vert _{B_{q,r,\gamma
}^{s}\left( R^{N};E\right) }+\sum\limits_{k=1}^{N}\left\Vert F^{-\shortmid
}\left( \xi _{k}^{l}\hat{u}\right) \right\Vert _{B_{q,r,\gamma }^{s}\left(
R^{N};E\right) }\right] .
\end{eqnarray*}%
\ Hence, the inequality (17) will be followed if we can prove the following
estimate \ 
\begin{equation*}
\left\Vert F^{-\shortmid }\left[ \left( i\xi \right) ^{\alpha }A^{1-x-\mu }%
\hat{u}\right] \right\Vert _{B_{q,r,\gamma }^{s}\left( R^{N};E\right) }\leq
C\left\Vert F^{-\shortmid }\left( \left[ A+I\theta \right] \hat{u}\right)
\right\Vert _{B_{q,r,\gamma }^{s}\left( R^{N};E\right) }\eqno(14)
\end{equation*}%
for all $u\in B_{q,r,\gamma }^{l,s}\left( R^{N};E\left( A\right) ,E\right) ,$
where 
\begin{equation*}
\theta =\theta \left( \xi \right) =\sum\limits_{k=1}^{N}\left\vert \xi
_{k}\right\vert ^{l}\in S\left( \varphi \right) .
\end{equation*}%
Let us express the left hand side of (18) as follows 
\begin{eqnarray*}
&&\left\Vert F^{-\shortmid }\left[ \left( i\xi \right) ^{\alpha }A^{1-x-\mu }%
\hat{u}\right] \right\Vert _{B_{q,r,\gamma }^{s}\left( R^{N};E\right) } \\
&=&\left\Vert F^{-\shortmid }\left( i\xi \right) ^{\alpha }A^{1-x-\mu }\left[
(A+I\theta \right] ^{-1}\left[ \left( A+I\theta \right) \right] \hat{u}%
\right\Vert _{B_{q,r,\gamma }^{s}\left( R^{N};E\right) }
\end{eqnarray*}%
(Since\ $A$ is the positive operator in\ $E$ and $\theta \left( \xi \right)
\in S\left( \varphi \right) ,$ $\left[ (A+I\theta \right] ^{-1}$ exists ).
From Corollary 2.18 we know that%
\begin{equation*}
\begin{array}{lll}
\left\Vert F^{-\shortmid }\left( i\xi \right) ^{\alpha }A^{1-x-\mu }\left[
(A+I\theta \right] ^{-1}\left[ \left( A+I\theta \right) \right] \hat{u}%
\right\Vert _{B_{q,r,\gamma }^{s}\left( R^{N};E\right) } & \leq &  \\ 
C\left\Vert F^{-\shortmid }\left( \left[ A+I\theta \right] \hat{u}\right)
\right\Vert _{B_{q,r,\gamma }^{s}\left( R^{N};E\right) } &  & 
\end{array}%
\end{equation*}%
holds if operator-function $\Psi \left( \xi \right) =\left( i\xi \right)
^{\alpha }A^{1-x-\mu }(A+\theta )^{-1}$ satisfies Mikhlin's condition for
each multi--index $\beta ,$ $|\beta |\leq \left\lceil \frac{N}{p}%
\right\rceil +1.$ It is clear that 
\begin{equation*}
\begin{array}{rll}
\left\Vert (1+|\xi |)^{|\beta |}D^{\beta }\Psi \left( \xi \right)
\right\Vert _{L_{\infty }(B(E))} & \leq & \displaystyle\sum\limits_{k=0}^{|%
\beta |}\left\Vert |\xi |^{k}D^{\beta }\Psi \left( \xi \right) \right\Vert
_{L_{\infty }(B(E))} \\ 
\vspace{-3mm} &  & 
\end{array}%
\end{equation*}%
Therefore, it is enough to show 
\begin{equation*}
\begin{array}{rll}
\left\Vert |\xi |^{k}D^{\beta }\Psi \left( \xi \right) \right\Vert
_{L_{\infty }(B(E))} & \leq & C \\ 
\vspace{-3mm} &  & 
\end{array}%
\end{equation*}%
for $k=0,1,\cdot \cdot \cdot |\beta |$ and $|\beta |\leq \left\lceil \frac{N%
}{p}\right\rceil +1.$ It is proven in [14] that $\Psi $ satisfies Miklin's
condition. Hence proof is completed.\hbox{\vrule
height7pt width5pt}

\section*{4. Degenerate differential--operator equations}

In this section we study degenerate elliptic DOE 
\begin{equation*}
(L+\lambda )u~=~-\left( \gamma (t)\frac{d}{dt}\right) ^{2}u+A_{1}(t)\left(
\gamma (t)\frac{d}{dt}\right) u+A_{\lambda }u~=~f\eqno(15)
\end{equation*}%
in $B_{q,r}^{s}(R;E),$ where $A_{\lambda }=A+\lambda I$ and $A_{1}(x)$ are
possible unbounded operators in a Banach space $E$. Let $E$ and $E_{0}$ be
Banach spaces such that $E_{0}$ is continuously and densely embedded in $E.$
Then 
\begin{equation*}
\begin{array}{l}
\displaystyle B_{p,q}^{\left[ l\right] ,s}(R;E_{0},E)~=~\left\{ u:u\in
B_{p,q}^{s}(R;E_{0}),~~D^{\left[ l\right] }u\in B_{p,q}^{s}(R;E)\right\} ,
\\ 
\vspace{-3mm} \\ 
\displaystyle\Vert u\Vert _{B_{p,q}^{\left[ l\right] ,s}(R;E_{0,}E)}~=~\Vert
u\Vert _{B_{p,q}^{s}(R;E_{0})}+\left\Vert D^{\left[ l\right] }u\right\Vert
_{B_{p,q}^{s}(R,E)}~<~\infty%
\end{array}%
\end{equation*}%
denotes the Besov-Lions spaces where 
\begin{equation*}
D^{[i]}~=~\left( \gamma (t)\frac{d}{dt}\right) ^{i}.
\end{equation*}

\vspace{3mm}

\textbf{Remark 4.1.} It is clear that under a substitution 
\begin{equation*}
\tau ~=~\int\limits_{0}^{t}\gamma ^{-1}(y)dy\eqno(16)
\end{equation*}%
spaces $B_{q,r}^{s}(R;E)$ and $B_{q,r}^{[2],s}(R;E(A),E),$ map
isomorphically onto the weighted spaces $B_{q,r,\tilde{\gamma}}^{s}(R;E)$ and%
$~B_{q,r,\tilde{\gamma}}^{2,s}(R;E(A),E)$ respectively, where $\tilde{\gamma}%
=\tilde{\gamma}(\tau )=\gamma (t(\tau )).$ Note that, (16) transforms
degenerate problem (15) in $B_{q,r}^{s}(R;E)$ to the following
non--degenerate problem

\begin{equation*}
(L+\lambda )u~=-u^{\prime \prime }+A_{1}(t)u^{\prime }+A_{\lambda }u~=f\eqno%
(17)
\end{equation*}%
in $B_{q,r,\tilde{\gamma}}^{s}(R;E).$

\vspace{3mm}

\textbf{Theorem 4.2.} Assume $\tilde{\gamma}$ satisfies the Condition 2. Let 
$E$ be a Banach space with Fourier type $p$, $A$ be a $\varphi $--positive
operator in $E$ for $\varphi \in \lbrack 0,\pi )$ and%
\begin{equation*}
A_{1}(\cdot )A^{-(\frac{1}{2}-\mu )}\in L_{\infty }(R,B(E)),\text{ }0<\mu <%
\frac{1}{2}.
\end{equation*}%
Then for all $f\in B_{q,r,\tilde{\gamma}}^{s}(R;E),~r,~q\in \lbrack 1,\infty
]$ and sufficiently large$~\lambda \in S(\varphi )$ (17) has a unique
solution $u\in B_{q,r,\tilde{\gamma}}^{2,s}(R;E(A),E)$ satisfying coercive
estimate 
\begin{equation*}
\left\Vert u^{\prime \prime }\right\Vert _{B_{q,r,\tilde{\gamma}%
}^{s}(R;E)}+\Vert A_{1}u\Vert _{B_{q,r,\tilde{\gamma}}^{s}(R;E)}+\Vert
Au\Vert _{B_{q,r,\tilde{\gamma}}^{s}(R;E)}\leq C\Vert f\Vert _{B_{q,r,\tilde{%
\gamma}}^{s}(R;E)}.\eqno(18)
\end{equation*}%
\textbf{Proof.} We first show maximal regularity result for the principal
part of (17) i.e.%
\begin{equation*}
(L_{0}+\lambda )u~=-u^{\prime \prime }+A_{\lambda }u~=f.
\end{equation*}%
\ By applying the Fourier transform, we obtain 
\begin{equation*}
(A+\xi ^{2}+\lambda )u^{\symbol{94}}(\xi )~=f^{\symbol{94}}(\xi ).
\end{equation*}%
Since $\xi ^{2}+\lambda \in S(\varphi )$ for all $\xi \in R$ and $A$ is a
positive operator, solutions are of the form 
\begin{equation*}
u(x)~=F^{-1}[A+\xi ^{2}+\lambda ]^{-1}f^{\symbol{94}}.\eqno(19)
\end{equation*}%
By using (19), we get 
\begin{equation*}
\begin{array}{lll}
\Vert Au\Vert _{B_{q,r,\tilde{\gamma}}^{s}(R;E)} & = & \displaystyle%
\left\Vert F^{-1}A(A+\xi ^{2}+\lambda )^{-1}f^{\symbol{94}}\right\Vert
_{B_{q,r,\tilde{\gamma}}^{s}(R;E)} \\ 
\vspace{-3mm} &  &  \\ 
\Vert u^{\prime \prime }\Vert _{B_{q,r,\tilde{\gamma}}^{s}(R;E)} & = & %
\displaystyle\left\Vert F^{-1}[\xi ^{2}(A+\xi ^{2}+\lambda )^{-1}f^{\symbol{%
94}}]\right\Vert _{B_{q,r,\tilde{\gamma}}^{s}(R;E)}.%
\end{array}%
\end{equation*}%
Therefore, it suffices to show that operator--function 
\begin{equation*}
\sigma (\xi )~=\xi ^{2}(A+\xi ^{2}+\lambda )^{-1},
\end{equation*}%
is uniformly bounded multiplier in $B_{q,r,\tilde{\gamma}}^{s}(R;E).$ Since $%
\sigma \in C^{2}(R,B(E))$ and 
\begin{equation*}
(1+|\xi |)^{k}D^{k}\sigma (\xi )\in L_{\infty }(R,B(E))
\end{equation*}%
for each $k=0,1,2,$ Corollary 2.18 guarantees us that $\sigma $ is a
uniformly bounded Fourier multiplier in $B_{q,r,\tilde{\gamma}}^{s}(R;E).$
Thus we obtain%
\begin{equation*}
\Vert L_{0}u\Vert _{B_{q,r,\tilde{\gamma}}^{s}(R;E)}\leq \left\Vert
u\right\Vert _{B_{q,r,\tilde{\gamma}}^{2,s}\left( R;E\left( A\right)
,E\right) }\leq C\Vert (L_{0}+\lambda )u\Vert _{B_{q,r,\tilde{\gamma}%
}^{s}(R;E)}.\eqno(20)
\end{equation*}%
The above estimate implies that $L_{0}+\lambda $ \ has a bounded inverse
acting from $B_{q,r,\tilde{\gamma}}^{s}(R;E)$ into $B_{q,r,\tilde{\gamma}%
}^{2,s}\left( R;E\left( A\right) ,E\right) .$ Next we try to estimate lower
order term $L_{1}u(t)=A_{1}\left( t\right) u^{\prime }(t).$ In fact, the
Theorem 3.1 ensures that for all$\ u\in B_{q,r,\tilde{\gamma}}^{2,s}\left(
R;E\left( A\right) ,E\right) ,$%
\begin{eqnarray*}
\left\Vert L_{1}u\right\Vert _{B_{q,r,\tilde{\gamma}}^{s}} &\leq &\left\Vert
A_{1}\left( t\right) A^{-(\frac{1}{2}-\mu )}\right\Vert _{L_{\infty
}}\left\Vert A^{\frac{1}{2}-\mu }u^{\prime }(t)\right\Vert _{B_{q,r,\tilde{%
\gamma}}^{s}} \\
&\leq &C_{0}C_{\mu }\left[ h^{\mu }\left\Vert u\right\Vert _{B_{q,r,\tilde{%
\gamma}}^{2,s}\left( R;E\left( A\right) ,E\right) }+h^{-\left( 1-\mu \right)
}\left\Vert u\right\Vert _{B_{q,r,\tilde{\gamma}}^{s}\left( R;E\right) }%
\right] .
\end{eqnarray*}%
It is also clear that 
\begin{eqnarray*}
\left\Vert u\right\Vert _{B_{q,r,\tilde{\gamma}}^{s}\left( R;E\right) } &=&%
\frac{1}{\lambda }\left\Vert \left( L_{0}+\lambda \right)
u-L_{0}u\right\Vert _{B_{q,r,\tilde{\gamma}}^{s}} \\
&\leq &\frac{1}{\lambda }\left[ \left\Vert \left( L_{0}+\lambda \right)
u\right\Vert _{B_{q,r,\tilde{\gamma}}^{s}}+\left\Vert u\right\Vert _{B_{q,r,%
\tilde{\gamma}}^{2,s}\left( R;E\left( A\right) ,E\right) }\right] ,
\end{eqnarray*}%
\ \ \ which in its turn implies the following estimate 
\begin{equation*}
\begin{array}{lll}
\vspace{-3mm}\left\Vert L_{1}u\right\Vert _{B_{q,r,\tilde{\gamma}}^{s}} & 
\leq & \displaystyle C_{1}\left[ (h^{\mu }+h^{-\left( 1-\mu \right) }\frac{1%
}{\lambda })\left\Vert u\right\Vert _{B_{q,r,\tilde{\gamma}}^{2,s}\left(
R;E\left( A\right) ,E\right) }\right. \\ 
&  &  \\ 
& + & \displaystyle\left. h^{-\left( 1-\mu \right) }\frac{1}{\lambda }%
\left\Vert \left( L_{0}+\lambda \right) u\right\Vert _{B_{q,r,\tilde{\gamma}%
}^{s}}\right] \\ 
&  &  \\ 
& \leq & \displaystyle C_{1}\left[ (h^{\mu }+h^{-\left( 1-\mu \right) }\frac{%
1}{\lambda })C\Vert (L_{0}+\lambda )u\Vert _{B_{q,r,\tilde{\gamma}%
}^{s}(R;E)}\right. \\ 
&  &  \\ 
& + & \left. h^{-\left( 1-\mu \right) }\frac{1}{\lambda }\left\Vert \left(
L_{0}+\lambda \right) u\right\Vert _{B_{q,r,\tilde{\gamma}}^{s}}\right] \leq
\Vert (L_{0}+\lambda )u\Vert _{B_{q,r,\tilde{\gamma}}^{s}(R;E)} \\ 
&  &  \\ 
\vspace{-3mm} & \times & \displaystyle\left[ CC_{1}h^{\mu
}+C_{1}(C+1)h^{-\left( 1-\mu \right) }\frac{1}{\lambda }\right] .%
\end{array}%
\eqno(21)
\end{equation*}%
~ \ \ Therefore choosing \ $h$ and $\lambda $ so that 
\begin{equation*}
CC_{1}h^{\mu }<1\text{ and }C_{1}(C+1)\left\vert \lambda \right\vert
^{-1}h^{-\left( 1-\mu \right) }<1
\end{equation*}%
we obtain 
\begin{equation*}
\ \ \left\Vert L_{1}\left( L_{0}+\lambda \right) ^{-1}\right\Vert _{B\left(
B_{q,r,\tilde{\gamma}}^{s}\left( R;E\right) \right) }<1.\eqno(22)
\end{equation*}%
In view of (20), (22) and the perturbation theory of linear operators, we
conclude that $L+\lambda =(L_{0}+\lambda )+L_{1}$ is invertible and its
inverse is continuous i.e. 
\begin{equation*}
(L+\lambda )^{-1}=\left( L_{0}+\lambda \right) ^{-1}\left[ I+L_{1}\left(
L_{0}+\lambda \right) ^{-1}\right] ^{-1}:B_{q,r,\tilde{\gamma}}^{s}\left(
R;E\right) \rightarrow B_{q,r,\tilde{\gamma}}^{2,s}\left( R;E\left( A\right)
,E\right) .
\end{equation*}%
Moreover, combining the estimates (20) and (21) we get (18). Hence proof is
completed. \hbox{\vrule height7pt width5pt}

\end{document}